\documentclass[a4paper,reqno,12pt]{amsart}
\usepackage{t1enc}
\usepackage[margin=1in]{geometry}
\usepackage{ifthen}
\usepackage{graphicx}

\newcommand{\R}{\mathbf{R}}

\newcommand{\pr}{\mathbf{P}}
\newcommand{\ex}{\mathbf{E}}

\usepackage{color}

\theoremstyle{plain}
\newtheorem{theorem}{Theorem}
\newtheorem{lemma}{Lemma}
\newtheorem{corollary}{Corollary}
\newtheorem{proposition}{Proposition}

\theoremstyle{definition}
\newtheorem{definition}{Definition}

\theoremstyle{remark}

\newcommand{\formula}[2][nolabel]
{\ifthenelse{\equal{#1}{nolabel}}
 {\begin{align*} #2 \end{align*}}
 {\ifthenelse{\equal{#1}{}}
  {\begin{align} #2 \end{align}}
  {\begin{align} \label{#1} #2 \end{align}}
 }
}

\renewcommand{\H}{\mathbb{H}}

\numberwithin{equation}{section}

\begin{document}

\title[Hitting hyperbolic half-space]{Hitting hyperbolic half-space}
\thanks{The authors were supported by MNiSW grant N~N201~373136}
\author{Jacek Ma{\l}ecki, Grzegorz Serafin}
\address{Jacek Ma{\l}ecki,  \\ Institute of Mathematics and Computer Science \\ Wroc{\l}aw University of Technology \\ ul. Wybrze{\.z}e Wyspia{\'n}\-skiego 27 \\ 50-370 Wroc{\l}aw, Poland}
\email{jacek.malecki@pwr.wroc.pl}
\address{Grzegorz Serafin \\ Institute of Mathematics and Computer Science \\ Wroc{\l}aw University of Technology \\ ul. Wybrze{\.z}e Wyspia{\'n}\-skiego 27 \\ 50-370 Wroc{\l}aw, Poland}
\email{grzegorz.serafin@pwr.wroc.pl}

\begin{abstract}
Let $X^{(\mu)}=\{X_t^{(\mu)};t\geq 0\}$, $\mu>0$, be the $n$-dimensional hyperbolic Brownian motion with drift, that is a diffusion on the real hyperbolic space $\H^n$ having the Laplace-Beltrami operator with drift as its generator. We prove the reflection principle for $X^{(\mu)}$ which enables us to study the process $X^{(\mu)}$ killed when exiting the hyperbolic half-space, that is the set $D=\{x\in\H^n:x_1>0\}$. We provide formulae, uniform estimates and describe asymptotic behavior of the Green function and the Poisson kernel of $D$ for the process $X^{(\mu)}$. Finally, we derive formula for the $\lambda$-Poisson kernel of the set $D$. 
\end{abstract}

\maketitle
\section{Introduction}
In recent years we have seen a considerable growth of interest in stochastic processes on hyperbolic spaces. In a series of papers \cite{BTF:2001}, \cite{BTFY:2001}, \cite{BGS:2007}, \cite{BM:2007}, \cite{Gruet:1996}, \cite{Matsumoto:2001}, \cite{Zak:2007} many important aspects of the hyperbolic Brownian motion were investigated, such as formulae and properties of transition density functions, harmonic measures and Green functions. The hyperbolic Brownian motion is a canonical diffusion on hyperbolic spaces. It has the half of the Laplace-Beltrami operator as its generator. Consequently, it is a hyperbolic analogue of classical Brownian motion on Euclidean space. From the other side it is closely related to geometric Brownian motion and Bessel processes, see \cite{BR:2006} and \cite{BMR3:2011}. Thus, it has some important applications in the risk theory in financial mathematics, see \cite{D:1990} and \cite{Yor:1992a}. 

The present paper is natural continuation of the research started in \cite{BGS:2007}, where the integral representation of the Poisson kernel for the interior of the horocycle, along with the resulting analysis of the asymptotic behavior was given. We consider the hyperbolic Brownian motion with drift $X^{(\mu)}$, $\mu>0$, in the half-space model $\H^n=\{x\in \R^n: x_n>0\}$ of $n$-dimensional real hyperbolic space, i.e. a diffusion with a generator $\frac{1}{2}\Delta_\mu$, where 
\formula{
\Delta_\mu  = x_n^2\sum_{k=1}^n\frac{\partial^2 }{\partial x_k^2}-(2\mu-1)x_n\frac{\partial }{\partial x_n}.
}
Note that for $\mu = (n-1)/2$ the process $X^{(\mu)}$ is a classical hyperbolic Brownian motion on $\H^n$.
We describe the Green function $G_D(x,y)$ and the density function of the harmonic measure (the Poisson kernel $P_D(x,y)$) of the hyperbolic half-space, i.e. the set $D=\{x\in\H^n:x_1>0\}$. We provide formulae for those objects as well as the uniform estimates for the full range of parameters $x$ and $y$. Note that the operator $\Delta_\mu$ is a strongly elliptic operator on every bounded subset of $\H^n$ (bounded in hyperbolic metric). Then, the general result of \cite{CranstonZhao:1987} implies comparability of the hyperbolic and Euclidean potential theory (comparability of corresponding Poisson kernels and Green functions). Since $D$ is unbounded, this kind of result is no longer true and deriving explicit formulae for $P_D$ and $G_D$ seems to be the only way to analyze the asymptotic behavior of those objects. It is also worth mentioning that the uniform estimates presented in Theorem \ref{thm:Green:estimate} and Theorem \ref{thm:Poisson:estimate} together with those obtained in \cite{BMR3:2011} are the only results of this kind in the case of unbounded subsets of hyperbolic spaces. 

The paper is organized as follows. In Preliminaries we provide basic information about Bessel processes and hyperbolic spaces. In the part related to Bessel processes we follow the exposition given in \cite{MatsumotoYor:2005a} and \cite{MatsumotoYor:2005b}, where we refer the Reader for more details and deeper insight into the subject (see also \cite{RevuzYor:2005}). In Section \ref{sec:HBM}, after describing the structure of the hyperbolic Brownian motion with drift $X^{(\mu)}$, we provide formulae for the transition density function, potentials and global Poisson kernel for the process $X^{(\mu)}$. In Proposition  \ref{prop:reflectionprinciple} we prove the reflection principle which gives the formulae for the transition density function for the process $X^{(\mu)}_D$, i.e. the process $X^{(\mu)}$ killed at the boundary of the hyperbolic half-space $D$. In Sections \ref{sec:GreenFunction} and \ref{sec:PoissonKernel} we introduce the formulae for the Green function and the Poisson kernel in terms of the modified Legendre functions and examine the asymptotic properties of those objects as well as its uniform estimates. All the properties of the Legendre functions used in those sections are given in the Appendix. We end the paper with Section \ref{sec:LambdaPK}, where an application of the Girsanov's Theorem together with the results of Section \ref{sec:PoissonKernel} lead to the formulae for the $\lambda$-Poisson kernel of the set $D$ and the solution of the appropriate Dirichlet problem on $D$ for the Laplace-Beltrami operator with drift.

\section{Preliminaries}
\label{sec:Preliminaries}
This section contains some preliminary material related to Bessel processes and geometric Brownian motion as well as real hyperbolic spaces which are indispensable in the sequel. 

\subsection{Bessel process and geometric Brownian motion}
We denote by $R^{(\nu)}=\{R^{(\nu)}_t,t\geq 0\}$ the Bessel process with index $\nu\in \R$ starting from $R_0^{(\nu)}=x>0$. In the case $-1<\nu<0$, it is when the point $0$ is non-singular, we impose killing condition at $0$. Then, for every $\nu<0$, the transition density function is given by (see \cite{BorodinSalminen:2002} p.134)
\formula[eq:bessel:density]{
g_t^{(\nu)}(w)=\frac x{t}\left(\frac{x}{w}\right)^\nu\exp\left(-\frac{x^2+w^2}{2t}\right)I_{-\nu}\left(\frac{xw}{t}\right)\quad x,w>0\/.
}
Here $I_\nu(z)$ is the modified Bessel function.

Let us denote by $B = \{B_t; t \geq 0\}$ the one-dimensional Brownian motion starting from $0$ and by $B^{(\nu)} = \{B^{(\nu)}_t = Bt+\nu t; t \geq 0\}$ the Brownian motion with a constant drift $\nu \in \R$. Let $Y_t^{(\nu)}$ be the geometric Brownian motion with drift starting from $x>0$, i.e.
\formula{
  Y^{(\nu)} = \{x\exp(B^{(\nu)}_t),t\geq 0\}\/.
}
Bessel process and geometric Brownian motion are related one to the other by the \textit{Lamperti relation}, which gives
\formula[eq:Lamperti]{
   Y^{(\nu)}_t \stackrel{d}{=}R^{(\nu)}_{A_x^{(\nu)}(t)}\/,
}
where the integral functional $A_x^{(\nu)}(t)$ is defined by
\formula[eq:At:defn]{
   A^{(\nu)}_{x}(t)=\int_0^t\left(Y^{(\nu)}(s)\right)^2ds=x^2\int_0^t\exp\left(2B_s+2\nu s\right)ds.
}
Whenever $\mu$ is strictly positive the limit of $A^{(-\mu)}_{x_n}(t)$ when $t$ goes to $\infty$ exists a.s. and its density function $h_{x_n}^{(\mu)}$ is given by (see \cite{D:1990})
\formula[eq:gestoscA]{
h_{x_n}^{(\mu)}(u)&=\frac{x_n^{2\mu}}{\Gamma(\mu)2^\mu}\frac{e^{-x_n^2/2u}}{u^{1+\mu}}\quad 1{\hskip -2.5 pt}\hbox{l}_{(0,\infty)}(u).
}

\subsection{Real hyperbolic space}
For every $n=1,2,3,\ldots$ we consider a half-space model of $n$-dimensional real hyperbolic space, i.e. the set $\H^n=\left\{x=(x_1,\ldots,x_n)\in\R^n:x_n>0\right\}$ with Riemannian metric
\formula{
ds^2 = \frac{dx^2}{x_n^2}.
}
The hyperbolic distance $d_{\H^n}(x,y)$ between two points $x$ and $y$ from $\H^n$ is described by
\formula[eq:coshdxy]{
\cosh(d_{\H^n}(x,y)) = 1+\frac{\left|x-y\right|^2}{2x_ny_n}.
}
The group of isometries of $\H^n$ is generated by reflexions with respect to hyperplanes $\{x\in\R^n: x_1=\textrm{const.}\}$ and inversions with respect to the spheres perpendicular to the hyperplane $P=\left\{x\in\R^n:x_n=0\right\}$. Thus, the set
\formula[eq:D:defn]{
D = \{x\in\H^n: x_1>0\}
}
can be considered as hyperbolic analogue of the half-space. Indeed, the set $D$ and the interior of the complement of $D$ can be transformed one to the other using the isometry $\H^n\ni x=(x_1,x_2,\ldots,x_n)\to (-x_1,x_2,\ldots,x_n)\in\H^n$. We will denote by $\partial_1D = \{x\in\H^n: x_1=0\}$ the boundary of $D$ in $\H^n$ and by $\partial_2D = \{x\in\R^n:x_n=0,x_1>0\}$ the second part of the boundary of $D$ considered as a subset of $\R^n$. Note that $\partial_2D$  does not belong to $\H^n$. We will write $\partial D$ for $\partial_1 D\cup \partial_2 D$.
The hyperbolic volume element is given by 
\formula[eq:volume]{
  dV_n= \frac{1}{x_n^n}dx_1...dx_n,
}
where $dx_1...dx_n$ stands for the Lebesgue measure on $\R^n$. Finally, the Laplace-Beltrami operator associated with the Riemannian metric is 
\formula[eq:LaplaceBeltrami]{
\Delta_{\H^n}  = x_n^2\sum_{k=1}^n\frac{\partial^2 }{\partial x_k^2}-(n-2)x_n\frac{\partial }{\partial x_n}.
}
This operator is the only second order elliptic differential operator on $\H^n$, annihilating constant functions, which is invariant under isometries of the space. 
\subsection{Notation}
For every $x=(x_1,\ldots,x_{n-1},x_n)\in\R^{n}$ we introduce the following notation
\formula{
\tilde{x}=(x_1,\ldots,x_{n-1})\in\R^{n-1}\/,\quad 
\bar{x} =(-x_1,x_2,\ldots,x_n)\in\R^{n} \/.}
Consequently, for every subset $A\subset\R^n$ we denote $\bar{A} = \{\bar{x}\in \R^n: x\in A\}$. Moreover, for every $x,y\in\H^n$ we put 
\formula{
\rho = d_{\H^n}(x,y)\/,\quad 
\bar{\rho} = d_{\H^n}(x,\bar{y})\/.
} 
Note that using (\ref{eq:coshdxy}) we easily obtain the following relation between $\rho$ and $\bar{\rho}$
\formula[eq:rhobarrho]{
\cosh\bar{\rho} = \cos \rho+\frac{2x_1y_1}{x_ny_n}\/,\quad x,y\in\H^n\/.
}
Finally, we write $f\stackrel{c}{\approx}g$, $x\in A$ whenever $c^{-1}f(x)\leq g(x)\leq c f(x)$ for every $x\in A$.


\section{Hyperbolic Brownian motion with drift and reflection principle}
\label{sec:HBM}
\subsection{Definition and structure}
For every $\mu > 0$, we define the following operator
\formula{
\Delta_\mu  := x_n^2\sum_{k=1}^n\frac{\partial^2 }{\partial x_k^2}-(2\mu-1)x_n\frac{\partial }{\partial x_n}.
}
Note that for $\mu=\frac{n-1}{2}$ the operator is the Laplace-Beltrami operator (\ref{eq:LaplaceBeltrami}).

The hyperbolic Brownian motion (HBM) with drift is a diffusion $X^{(\mu)} = \{X^{(\mu)}_t;t\geq0 \}$ on $\H^n$ starting from $X_0^{(\mu)}=x\in\H^n$
with a generator $\frac12 \Delta_\mu$. For $\mu=\frac{n-1}{2}$ we obtain the standard HBM on $\H^n$. 

The structure of $X^{(\mu)}$ can be described as follows. 
If $\beta=\{(\beta_1(t),\ldots,\beta_{n-1}(t)),t\geq 0\}$ is $(n-1)$-dimensional Brownian motion starting from $(x_1,\ldots,x_{n-1})\in\R^{n-1}$ independent from a geometric Brownian motion $Y^{(-\mu)}$ starting from $x_n>0$ then 
\formula[eq:HBM:GBM]{
   X^{(\mu)}_t \stackrel{d}{=} (\beta_1(A_{x_n}^{(-\mu)}(t)),\ldots,\beta_{n-1}(A_{x_n}^{(-\mu)}(t)),Y^{(-\mu)}_t) 
}
Moreover, following the \textit{Lamperti relation} one can also show that 
\formula[eq:HBM:bessel]{
  X^{(\mu)}_t\stackrel{d}{=}\left(\beta_1(A^{(-\mu)}_{x_n}(t)),...,\beta_{n-1}(A^{(-\mu)}_{x_n}(t)),R^{(-\mu)}(A^{(-\mu)}_{x_n}(t))\right),
}
where $R^{(-\mu)}=\{R^{(-\mu)}_t,t\geq t\}$ is a Bessel motion with index $-\mu$ starting from $x_n>0$ independent from Brownian motion $\beta$.

\subsection{Transition density function, potentials and global Poisson kernel}
In this section we compute basic characteristics of the process $X^{(\mu)}$ such as a transition density function, its potential kernel and so-called global Poisson kernel, i.e. the density function of $\lim_{t\to \infty }X_t^{(\mu)}$.
\begin{proposition}
\label{prop:pt:formula}
Transition density function $p^{(\mu)}(t,x,y)$ of $X^{(\mu)}$ with respect to the canonical volume element $dV_n$ is 
\formula[eq:pt:formula]{
p^{(\mu)}(t,x,y)=\left(\frac{x_n}{y_n}\right)^{\mu-\nu}\frac{e^{-\mu^2t/2}}{\pi(2\pi)^{n/2}t^{1/2}}\Gamma\left(\frac{n+1}{2}\right)\int_0^\infty\frac{e^{(\pi^2-u^2)/2t}\sinh u\sin(\pi u/t)}{\left(\cosh u+\cosh \rho\right)^{(n+1)/2}}du,
}
where $\nu=\frac{n-1}{2}$ and $\rho=d_{\H^n}(x,y)$.
\end{proposition}
\begin{proof}
The density function $f_{x_n,t}^{(\mu)}$ of $(A^{(-\mu)}_{x_n}(t),Y_t^{(-\mu)})$ was computed in \cite{Yor:1992a}. We have
\formula[eq:densityf]{
f_{x_n,t}^{(\mu)}(u,v)=\left(\frac{v}{x_n}\right)^{-\mu}e^{-\mu^2t/2}\frac1{uv} \exp\left(-\frac{x_n^2+v^2}{2u}\right)\theta_{x_nv/u}(t),
}
where
\formula[eq:theta]{
\theta_r(t)=\frac{r}{(2\pi^3t)^{1/2}}\int_0^\infty e^{(\pi^2-b^2)/(2t)}e^{-r\cosh (b)}\sinh(b)\sin\left(\frac{\pi b}{t}\right)db.
}
Moreover, the Laplace transform of the function $\theta_r$ is (see \cite{Yor:1980})
\formula[eq:theta:Laplace]{
\int_0^\infty e^{-\lambda t}\theta_r(t)dt=I_{\sqrt{2\lambda}}(r),
}
where $I_\nu(z)$ is a modified Bessel function.
Relation (\ref{eq:HBM:GBM}) and the fact that Brownian motion $\beta$ and geometric Brownian motion $Y^{(-\mu)}$ are independent allow us to write
\formula{
p^{(\mu)}(t,x,y)& =y_n^n\int_0^\infty \frac{1}{(2\pi u)^{(n-1)/2}}\exp\left(-\frac1{2u}\sum_{k=1}^{n-1}(x_k-y_k)^2\right)f_{x_n,t}^{(\mu)}(u,y_n)du\/.
}
The factor $y_n^n$ appearing in front of the integral follows from formula (\ref{eq:volume}). Consequently, using the formulae (\ref{eq:densityf}), (\ref{eq:coshdxy}) and substituting $u={x_ny_n}/{s}$ we obtain 
\formula[eq:gestosc:theta]{
p^{(\mu)}(t,x,y)&=\left(\frac{x_n}{y_n}\right)^{\mu-\nu}\frac{e^{-\mu^2t/2}}{(2\pi)^{\nu}}\int_0^\infty s^{\frac{n-3}{2}} e^{-s\cosh\rho}\theta_{s}(t)ds\/.
}
The formula (\ref{eq:theta}) and Fubini theorem lead to 
\formula{
\int_0^\infty s^{\frac{n-3}{2}} e^{-s\cosh\rho}\theta_{s}(t)ds = \int_0^\infty \frac{e^{(\pi^2-b^2)/(2t)}}{(2\pi^3t)^{1/2}}\sinh(b)\sin\left(\frac{\pi b}{t}\right)\int_0^\infty s^{\frac{n-1}{2}} e^{-s(\cosh\rho+\cosh b)}dsdb
}
and the substitution $s=\frac{w}{\cosh\rho+\cosh b}$ in the inner integral gives
\formula{
p^{(\mu)}(t,x,y)&=\left(\frac{x_n}{y_n}\right)^{\mu-\nu}\frac{e^{-\mu^2t/2}}{\pi(2\pi)^{n/2}t^{1/2}}\int_0^\infty\frac{e^{(\pi^2-b^2)/2t}\sinh b\sin(\pi b/t)}{\left(\cosh b+\cosh \rho\right)^{(n+1)/2}}\int_0^\infty w^{\frac{n-1}{2}} e^{-w}dw\,db\\
&=\left(\frac{x_n}{y_n}\right)^{\mu-\nu}\frac{e^{-\mu^2t/2}}{\pi(2\pi)^{n/2}t^{1/2}}\Gamma\left(\frac{n+1}{2}\right)\int_0^\infty\frac{e^{(\pi^2-b^2)/2t}\sinh b\sin(\pi b/t)}{\left(\cosh b+\cosh \rho\right)^{(n+1)/2}}db.
}
\end{proof}
An easy consequence of Proposition \ref{prop:pt:formula} is the following relation between transition density functions for the processes with different indices.
\begin{corollary}For every $x,y\in \H^n$ we have
\begin{equation}
\label{zaleznoscgestosci}
p^{(\mu)}(t,x,y)=\left(\frac{x_n}{y_n}\right)^{\mu-\nu}\exp\left(-\frac{\mu^2-\nu^2}{2}\,t\right)p^{\left(\frac{n-1}{2}\right)}(t,x,y)\/,
\end{equation}
where $\nu=\frac{n-1}{2}$.
\end{corollary}
Another consequence of the formula (\ref{eq:pt:formula}) is the representation of the potential kernel of the process $X^{(\mu)}$, which is defined by
\formula{
V^{(\mu)}(x,y) = \int_0^\infty p^{(\mu)}(t,x,y)dt\/,\quad x,y\in\H^n\/,
}
in terms of modified Legendre functions.
\begin{proposition}
For every $x,y\in\H^n$ we have
\formula[eq:potential:formula]{
V^{(\mu)}(x,y)=\left(\frac{x_n}{y_n}\right)^{\mu-\nu}\frac{2e^{-(\nu-1/2)i\pi }}{(2\pi)^{n/2}} \frac{Q^{\nu-1/2}_{\mu-1/2}(\cosh \rho)}{\sinh^{\nu-1/2}\rho},
}
where $\nu=\frac{n-1}{2}$ and $\rho=d_{\H^n}(x,y)$.
\end{proposition}
\begin{proof}
Representation (\ref{eq:gestosc:theta}) of $p^{(\mu)}(t,x,y)$ in terms of the function $\theta_r$ together with the formula (\ref{eq:theta:Laplace}) for its Laplace transform imply
\formula{
\int_0^\infty p^{(\mu)}(t,x,y)dt &{=} \left(\frac{x_n}{y_n}\right)^{\mu-\nu}\frac{1}{(2\pi)^{\nu}}\int_0^\infty s^{\frac{n-3}{2}} e^{-s\cosh\rho}\int_0^\infty e^{-\mu^2t/2}\theta_{s}(t)dt\;ds\\
&{=}\left(\frac{x_n}{y_n}\right)^{\mu-\nu}\frac{1}{(2\pi)^{\nu}}\int_0^\infty s^{\nu-1/2} e^{-s\cosh\rho}I_\mu(s)ds.
}
The integral formula (see \cite{GradsteinRyzhik:2007} 6.622 p. 702) 
\formula[eq:gr6.622]{
\int_0^\infty s^{\nu-1} e^{-s\cosh\rho}I_\mu(s)ds=\sqrt{\frac2{\pi}}e^{-(\nu-1/2)i\pi } \frac{Q^{\nu-1/2}_{\mu-1/2}(\cosh \rho)}{\sinh^{\nu-1/2}\rho},
}
ends the proof.

\end{proof}
Finally, the representation of the Legendre functions in terms of elementary function provided in Appendix gives
\begin{corollary} For odd $n\geq3$ we have
\formula{
V^{(\mu)}(x,y) =\left(\frac{x_n}{y_n}\right)^{\mu-\nu}\frac{\Gamma(\mu+\nu)}{(2\pi)^{\nu}} \sum_{k=0}^{\frac{n-3}2}\frac{a(\frac{n-3}2,k)}{\Gamma(\mu+1+k)}\frac{(\cosh\rho-\sinh\rho)^{\mu+k}}{(\sinh\rho)^{1/2+k}},
}
where $\nu=\frac{n-1}{2}$, $\rho=d_{\H^n}(x,y)$ and $a(m,k)$ are the coefficients of the Bessel polynomials
\formula{
 a(m,k)=\frac{(m+k)!}{(m-k)!k!2^k}.
 }
\end{corollary}
\begin{proof}
Follows immediately from (\ref{suma}) and (\ref{eq:potential:formula}).
\end{proof}

We denote by $P^{(\mu)}(x,y)$ the global Poisson kernel, i.e. the density function of the random variable $X^{(\mu)}_\infty=\lim_{t\to\infty}X^{(\mu)}_t$, where $x=X^{(\mu)}_0$ and $y\in\R^{n-1}\times\left\{0\right\}$. The existence of this limit follows from the representation (\ref{eq:HBM:GBM}) and the fact that $Y^{(-\mu)}_t$ tends to zero a.s. whenever $\mu>0$.
\begin{proposition}
\label{prop:gpk:formula} 
For $x\in\H^n$, $y\in\R^{n-1}\times\left\{0\right\}$ we have
\begin{equation}
\label{eq:gpk:formula}
P^{(\mu)}(x,y)=\frac{\Gamma(\mu+\nu)}{\Gamma(\mu)\pi^\nu }\frac{x_n^{2\mu}}{\left| x- y\right|^{2\mu+2\nu}},
\end{equation}
where $\nu=\frac{n-1}{2}$.
\end{proposition}

\begin{proof}
Using the representation (\ref{eq:HBM:GBM}), the fact that $Y^{(-\mu)}_\infty =0$ a.s and the formula (\ref{eq:gestoscA}) we obtain
\formula{
P^{(\mu)}(x,y)&=\int_{0}^\infty \frac{1}{(\sqrt{2\pi u})^{n-1}}\exp\left(-{\frac{\left|\tilde x-\tilde y\right|^2}{2u}}\right)h_{x_n}^{(\mu)}(u)du\\
&=\frac{x_n^{2\mu}}{\Gamma(\mu)2^\mu(\sqrt{2\pi})^{n-1}}\int_{0}^\infty \exp\left(-{\frac{\left|\tilde x-\tilde y\right|^2+x_n^2}{2u}}\right){u^{-(1+\mu+\nu)}}du.
}
Making the substitution $u=\frac{\left|\tilde x-\tilde y\right|+x_n^2}{2v}$ in the last integral we arrive at
\formula{P^{(\mu)}(x,y)=&\frac{x_n^{2\mu}}{\Gamma(\mu)\pi^\nu \left(\left|\tilde x-\tilde y\right|^2+x_n^2\right)^{\mu+\nu}}\int_{0}^\infty e^{-v}v^{\mu+\nu-1}dv
=\frac{\Gamma(\mu+\nu)}{\Gamma(\mu)\pi^\nu }\frac{x_n^{2\mu}}{\left| x- y\right|^{2\mu+2\nu}}.
}
\end{proof}

\subsection{Reflection principle}
We denote by $\tau_D = \inf\{t>0:X^{(\mu)}_t\notin D\}$ the first exit time of the process $X^{(\mu)}$ from the hyperbolic half-space $D$ (see (\ref{eq:D:defn}) for its definition) and $X^{(\mu)}_D$ stands for the process $X^{(\mu)}$ killed when exiting $D$. Its transition density function is given by the Hunt formula
\formula[eq:ptd:defn]{
p_D^{(\mu)}(t,x,y) = p^{(\mu)}(t,x,y)-\ex^x[t<\tau_D;p^{(\mu)}(t-\tau_D,X^{(\mu)},y)]\/,\quad t>0,\ \ x,y\in D\/.
}
We introduce the following reflexion principle for HBM with drift. Recall that using the notation introduced in the Preliminaries we have $\overline{X^{(\mu)}}= (-X_1^{(\mu)},X_2^{(\mu)},\ldots,X_n^{(\mu)})$.
\begin{proposition}[Reflection principle]
\label{prop:reflectionprinciple}
Let us define
\begin{equation}
  Z^{(\mu)}_t = 
    \left\{
      \begin{array}{rl}
         X^{(\mu)}_t,&  \tau_D>t \\
         \overline{X^{(\mu)}_t},& \tau_D\leq t . \\
       \end{array} \right.
\end{equation}
then $\{Z^{(\mu)}_t\}_{t\geq0}\stackrel{d}{=}\{X^{(\mu)}_t\}_{t\geq0}$.

\end{proposition}

\begin{proof} 
Fix $m=1,2,\ldots$,  $0=t_0<t_1<\ldots<t_m<t_{m+1}=\infty$ and $A_1,\ldots, A_m \in\mathcal{B}(\H^n)$. Then, for every $k=0,1,\ldots,m$ we have
\formula{
&\pr^x(t_k\leq \tau_D<t_{k+1}, Z^{(\mu)}_{t_1}\in A_1,\ldots, Z^{(\mu)}_{t_m}\in A_m)= \\
&= \pr^x(t_k\leq \tau_D<t_{k+1}, Z^{(\mu)}_{t_1}\in A_1,\ldots, Z^{(\mu)}_{t_k}\in A_k, Z^{(\mu)}_{t_{k+1}}\in A_{k+1},\ldots,Z^{(\mu)}_{t_m}\in A_m)\\
&= \pr^x(t_k\leq \tau_D<t_{k+1}, X^{(\mu)}_{t_1}\in A_1,\ldots, X^{(\mu)}_{t_k}\in A_k, X^{(\mu)}_{t_{k+1}}\in \overline{A_{k+1}},\ldots,X^{(\mu)}_{t_m}\in \overline{A_m})
}
Applying strong Markov property we obtain that the above-given expression is equal to 
\formula{
\ex^x(t_k\leq \tau_D<t_{k+1}, X^{(\mu)}_{t_1}\in A_1,\ldots, X^{(\mu)}_{t_k}\in A_k,\ex^{X^{(\mu)}_{\tau_D}}(X^{(\mu)}_{t_{k+1}-\tau_D}\in \overline{A_{k+1}},\ldots,X^{(\mu)}_{t_m-\tau_D}\in \overline{A_m}))
}
Note that the first coordinate of $X^{(\mu)}_{\tau_D}$ is equal to zero a.s. and for every $t>0$ we have $\pr^{(0,x_2,\ldots,x_n)}(X^{(\mu)}_t\in \overline{A}) = \pr^{(0,x_2,\ldots,x_n)}(X^{(\mu)}_t\in {A})$. Thus, we can change $\overline{A_{k+1}},\ldots,\overline{A_{m}}$ into ${A}_{k+1},\ldots,{A}_{m}$ and usage of the strong Markov property leads to 
\formula{
&\pr^x(t_k\leq \tau_D<t_{k+1}, Z^{(\mu)}_{t_1}\in A_1,\ldots, Z^{(\mu)}_{t_m}\in A_m) \\
&= \pr^x(t_k\leq \tau_D<t_{k+1}, X^{(\mu)}_{t_1}\in A_1,\ldots, X^{(\mu)}_{t_m}\in A_m)
}
Summing up the above-given equality with respect to $k$ we obtain that the finite dimensional distributions of the processes $X^{(\mu)}$ and $Z^{(\mu)}$ are the same and it ends the proof.
\end{proof}
The consequence of the reflection principle, as in the real case, is the following formula for transition density function of the killed process.
\begin{theorem}
\label{thm:pdt:formula}
For every $x,y\in D$ and $t> 0$ we have 
\formula[eq:pdt:formula]{
p^{(\mu)}_D(t,x,y) = p^{(\mu)}(t,x,y)-p^{(\mu)}(t,x,\bar{y})\/.}
\end{theorem}
\begin{proof}
Whenever $A\in\mathcal{B}(D)$ and $t>0$ we can write
\formula{
\pr^x(Z^{(\mu)}_t\in A)&= \pr^x(t<\tau_D;X^{(\mu)}_t\in A)+ \pr^x(t\geq \tau_D; \overline{X^{(\mu)}_t}\in A)\\
&= \pr^x(t<\tau_D;X^{(\mu)}_t\in A)+ \pr^x(\overline{X^{(\mu)}_t}\in A)\\
&= \pr^x(t<\tau_D;X^{(\mu)}_t\in A)+ \pr^x({X^{(\mu)}_t}\in \overline{A}).
}
Consequently, applying Proposition \ref{prop:reflectionprinciple}, we obtain
\formula{
\int_Ap_D^{(\mu)}(t,x,y)dV_n(y)&=\int_Ap_D^{(\mu)}(t,x,y)dV_n(y)-\int_{\bar A}p_D^{(\mu)}(t,x,y)dV_n(y)\\
&=\int_Ap_D^{(\mu)}(t,x,y)dV_n(y)-\int_Ap_D^{(\mu)}(t,x,\bar y)dV_n(y).
}
\end{proof}


\section{Green function}
\label{sec:GreenFunction}
The Green function of the set $D$ is defined as usual by
\formula{
G_D^{(\mu)}(x,y) = \int_0^\infty p_D^{(\mu)}(t,x,y)\/,\quad x,y \in D\/.
}
As an immediate consequence of Theorem \ref{thm:pdt:formula} we get the following representation for $G_D^{(\mu)}(x,y)$ in terms of the potential kernel $V^{(\mu)}(x,y)$.

\begin{theorem}
For every $x,y \in D$ we have
\label{thm:gf:formula}
\formula[eq:gf:formula]{
G_D^{(\mu)}(x,y) = V^{(\mu)}(x,y)-V^{(\mu)}(x,\bar{y})\/.
}
\end{theorem}
In the Propositions \ref{prop:gf:partial1}-\ref{prop:asym:last} we examine the asymptotic behavior of $G^{(\mu)}_D(x,y)$ at the boundary of the set $D$ and whenever $x-y$ tends to $0$.
\begin{proposition}
\label{prop:gf:partial1}
Let $z\in\partial_1D$, $x\in D$. Then
\formula[eq:gf:partial1]{
\lim_{\substack{y\rightarrow z\\y\in D}}y_1^{-1}G_D^{(\mu)}(x,y)= \frac{4e^{-in\pi/2 }}{(2\pi)^{n/2}} \frac{x_1}{x_nz_n}\left(\frac{x_n}{z_n}\right)^{\mu-\nu}\frac{Q^{\nu+1/2}_{\mu-1/2}(\cosh d_{\H^n}(x,z))}{(\sinh d_{\H^n}(x,z))^{\nu+1/2}}.
}
\end{proposition}

\begin{proof}
Using the representations (\ref{eq:gf:formula}), (\ref{eq:potential:formula}), the mean value  theorem and the formula (\ref{eq:legendre:diff}) we obtain 
\formula{
G_D^{(\mu)}(x,y)&=V^{(\mu)}(x,y)-V^{(\mu)}(x,\bar{y})\\
&=-\frac{4e^{-(\nu-1/2)i\pi }}{(2\pi)^{n/2}}\frac{x_1y_1}{x_ny_n} \left(\frac{x_n}{y_n}\right)^{\mu-\nu}\frac{\dfrac{Q^{\nu-1/2}_{\mu-1/2}(\cosh\bar{\rho})}{(\sinh \bar{\rho})^{\nu-1/2}}-\dfrac{Q^{\nu-1/2}_{\mu-1/2}(\cosh\rho)}{(\sinh \rho)^{\nu-1/2}}}{\cosh\bar{\rho}-\cosh\rho}\\
&=\frac{4e^{-in\pi/2 }}{(2\pi)^{n/2}} \frac{x_1y_1}{x_ny_n}\left(\frac{x_n}{y_n}\right)^{\mu-\nu}\frac{Q^{\nu+1/2}_{\mu-1/2}(\theta)}{(\sqrt{\theta^2-1})^{n/2}},
}
where $\cosh\rho\leq\theta\leq\cosh\bar{\rho}$. Using the fact that the function $z\to e^{-in\pi/2} Q^{\nu+1/2}_{\mu-1/2}(z)$ is decreasing on $(1,\infty)$ we get
\formula[eq:gf:inequality]{
\frac{4e^{-in\pi/2 }}{(2\pi)^{n/2}}\frac{Q^{\nu+1/2}_{\mu-1/2}(\cosh\bar{\rho})}{(\sinh\bar{\rho})^{n/2}}\leq \left(\frac{x_n}{y_n}\right)^{\nu-\mu}\frac{x_ny_n}{x_1y_1}\,G_D^{(\mu)}(x,y)\leq\frac{4e^{-in\pi/2 }}{(2\pi)^{n/2}} \frac{Q^{\nu+1/2}_{\mu-1/2}(\cosh\rho)}{(\sinh\rho)^{n/2}}.
}
By (\ref{eq:rhobarrho}) we get $\lim_{y\to z,y\in D}\cosh\bar{\rho}=\cosh\rho$ and consequently, by continuity of the Legendre function, we obtain (\ref{eq:gf:partial1}).
\end{proof}

\begin{proposition}
\label{prop:gf:partial2}
Let $z\in\partial_2D$, $x\in D$. Then 
\formula[eq:gf:partial2]{
\lim_{\substack{y\rightarrow z\\y\in D}}y_n^{1-n}G_D^{(\mu)}(x,y)=\frac{\Gamma(\mu+\nu)}{\pi^{\nu}\Gamma(\mu+1)}\left(\frac{x_n^{2\mu}}{\left|x-z\right|^{2\mu+2\nu}}-\frac{x_n^{2\mu}}{\left|x-\bar z\right|^{2\mu+2\nu}}\right).
}
\end{proposition}
\begin{proof}
Using (\ref{eq:potential:formula}) and (\ref{eq:coshdxy}) we get
\formula{
y_n^{1-n}V^{(\mu)}(x,y) &= x_n^{\mu-\nu}\frac{2e^{-(\nu-1/2)i\pi }}{(2\pi)^{n/2}}\,y_n^{-\mu-\nu}\frac{Q^{\nu-1/2}_{\mu-1/2}(\cosh \rho)}{\sinh^{\nu-1/2}\rho}\\
& = x_n^{\mu-\nu}\frac{2e^{-(\nu-1/2)i\pi }}{(2\pi)^{n/2}}\left(\frac{1}{y_n+\frac{\left|x-y\right|^2}{2x_n}}\right)^{\mu+\nu}\left(\cosh\rho^{\mu+\nu}\frac{Q^{\nu-1/2}_{\mu-1/2}(\cosh \rho)}{\sinh^{\nu-1/2}\rho}\right)
}
The asymptotic formula (\ref{eq:legendre:asymp04}) gives
\formula{
\lim_{\substack{y\rightarrow z\\y\in D}}y_n^{1-n}V^{(\mu)}(x,y) = \frac{\Gamma(\mu+\nu)}{\pi^{\nu}\Gamma(\mu+1)}\frac{x_n^{2\mu}}{\left|x-z\right|^{2\mu+2\nu}}.
}
Changing $y$ to $\bar{y}$ in above given equality and using (\ref{eq:gf:formula}) we prove (\ref{eq:gf:partial2}).
\end{proof}

\begin{proposition}
For $x\in D$, $(y_1,...,y_{n-1})\in\R^{n-1}$ we have
\formula[eq:gf:yninfty]{
\lim_{y_n\rightarrow\infty}y_n^{2\mu+2}G_D^{(\mu)}(x,y)=\frac{4\Gamma(\mu+\nu+1)}{\pi^\nu\Gamma(\mu+1)} x_n^{2\mu}x_1y_1.
}
\end{proposition}

\begin{proof}
Multiplying (\ref{eq:gf:inequality}) by $y_n^{\mu+\nu+1}$ we obtain
\formula[eq:gf:inequality2]{
\frac{4y_n^{\mu+\nu+1}}{(2\pi)^{n/2}} \frac{Q^{\nu+1/2}_{\mu-1/2}(\cosh\bar{\rho})}{(\sinh\bar{\rho})^{n/2}}\leq \frac{x_n^{\nu-\mu+1}y_n^{2\mu+2}}{x_1y_1e^{-n\pi i/2}}\,G_D^{(\mu)}(x,y)\leq\frac{4y_n^{\mu+\nu+1}}{(2\pi)^{n/2}} \frac{Q^{\nu+1/2}_{\mu-1/2}(\cosh\rho)}{(\sinh\rho)^{n/2}}.
}
By formula for hyperbolic distance (\ref{eq:coshdxy}) we have 
\formula[eq:gf:coshrho/yn]{
\lim_{y_n\rightarrow\infty}\frac{y_n}{\cosh{\rho}}=\lim_{y_n\rightarrow\infty}\left(\frac{1}{y_n}+\frac{1}{2x_n}\frac{\left|x-y\right|^2}{y_n^2}\right)^{-1}=2x_n.
}
Rewriting the right-hand side of (\ref{eq:gf:inequality2}) in the form
\formula{ \frac{4}{(2\pi)^{n/2}}\left(\frac{y_n}{\cosh{\rho}}\right)^{\mu+\nu+1}\left((\cosh{\rho})^{\mu+\nu+1}\frac{Q^{\nu+1/2}_{\mu-1/2}(\cosh{\rho})}{(\sinh{\rho})^{\nu+1/2}}\right)
}
and applying (\ref{eq:legendre:asymp04}) and (\ref{eq:gf:coshrho/yn}) we get that it tends to 
\formula{
\frac{4\Gamma(\mu+\nu+1)}{\pi^{\nu}\Gamma(\mu+1)} e^{-n\pi i/2} x_n^{\mu+\nu+1}
}
as $y_n\to \infty$. By (\ref{eq:rhobarrho}) we obtain $\lim_{y_n\rightarrow\infty}{\cosh\rho}/{\cosh\bar{\rho}} = 1$ and the limit of the left-hand side of (\ref{eq:gf:inequality2}) is exactly the same as above and consequently (\ref{eq:gf:yninfty}) holds.
\end{proof}

\begin{proposition}
\label{prop:asym:last}
For fixed $x\in D$ we have
\formula{
\lim_{y\rightarrow x}\left|x-y\right|^{n-2}G_D^{(\mu)}(x,y)=\frac{x_n^{n-2}}{2\pi^{n/2}}\Gamma\left(\frac{n-2}{2}\right).
}
\end{proposition}

\begin{proof}
Whenever $x\in D$ (i.e. $x_1\neq 0$) we have $\bar{x}\neq x$ and $\lim_{y\to x}V^{(\mu)}(x,\bar{y}) = V^{(\mu)}(x,\bar{x})$ which is finite.
Consequently
\formula{
\lim_{y\to x}\frac{G^{(\mu)}_D(x,y)}{|x-y|^{2-n}} &= \lim_{y\to x}\frac{V^{(\mu)}(x,y)}{|x-y|^{2-n}}\\
& = \frac{2e^{-in\pi/2 }}{(2\pi)^{n/2}} \lim_{y\rightarrow x}(2x_ny_n)^{\frac{n-2}{2}} \left(\frac{x_n}{y_n}\right)^{\mu-\nu}\left(\frac{\left|x-y\right|^2}{2x_ny_n}\right)^{\frac{n-2}{2}}\frac{Q^{\nu-1/2}_{\mu-1/2}(\cosh \rho)}{(\sinh\rho)^{\nu-1/2}}\\
&=\frac{^{-in\pi/2 }}{\pi^{n/2}}\lim_{y\rightarrow x}(x_ny_n)^{\frac{n-2}{2}}\left(\frac{x_n}{y_n}\right)^{\mu-\nu} \left(\cosh\rho-1\right)^{\nu-1/2}\frac{Q^{\nu-1/2}_{\mu-1/2}(\cosh \rho)}{(\sinh\rho)^{\nu-1/2}}\\
&{=}\frac{x_n^{n-2}}{2\pi^{n/2}}\Gamma\left(\frac{n-2}{2}\right),
}
where the last equality follows from (\ref{eq:legendre:asymp03}).
\end{proof}
We end this section by the following uniform estimates of the Green function $G_D^{(\mu)}(x,y)$. 
\begin{theorem}
\label{thm:Green:estimate}
There exists constant $c=c(n,\mu)>0$ such that 
\formula{
G_D^{(\mu)}(x,y)\stackrel{c}{\approx}\left(\frac{x_n}{y_n}\right)^{\mu-\nu}\frac{x_1y_1}{\left|x-\bar y\right|^2}\frac{(x_ny_n)^{\nu-1/2}}{\left|x-y\right|^{n-2}(\cosh \rho)^{\mu+1/2}},\quad x,y \in D\/,
}
where $\rho=d_{\H^n}(x,y)$, $\nu=\frac{n-1}2$.
\end{theorem}
\begin{proof}Using formulas (\ref{eq:gf:formula}), (\ref{eq:potential:formula}) and applying Lemma \ref{q-q} we get
\begin{eqnarray*}
G_D^{(\mu)}(x,y) &=& \left(\frac{x_n}{y_n}\right)^{\mu-\nu}\frac{2e^{-(\nu-1/2)i\pi }}{(2\pi)^{n/2}}\left[ \frac{Q^{\nu-1/2}_{\mu-1/2}(\cosh \rho)}{\sinh^{\nu-1/2}\rho}-\frac{Q^{\nu-1/2}_{\mu-1/2}(\cosh \bar\rho)}{\sinh^{\nu-1/2}\bar\rho}\right]\\
&\approx&\left(\frac{x_n}{y_n}\right)^{\mu-\nu}\frac{\cosh \bar\rho-\cosh \rho}{\cosh \bar\rho-1}\frac1{(\cosh \rho-1)^{\nu-1/2}(\cosh \rho)^{\mu+1/2}}\\
&=&2^{\nu+3/2}\left(\frac{x_n}{y_n}\right)^{\mu-\nu}\frac{x_1y_1}{\left|x-\bar y\right|^2}\frac{(x_ny_n)^{\nu-1/2}}{\left|x-y\right|^{n-2}(\cosh \rho)^{\mu+1/2}}.\\
\end{eqnarray*}
\end{proof}

\section{Poisson kernel}
\label{sec:PoissonKernel}
The random variable $X^{(\mu)}_{\tau_D}$ is supported on $\partial_1 D$ whenever $\tau_D<\infty$ and if $\tau_D=\infty$ then $X^{(\mu)}_{\tau_D}\in \partial_2 D$. We denote by $P^{(\mu)}_D(x,y)$ the Poisson kernel of $D$, i.e. the density function of random variable $X^{(\mu)}_{\tau_D}$ with respect to the Lebesgue measure (on $\partial_1 D$ and $\partial_2 D$). In the following theorem we provide explicit formula for $P^{(\mu)}_D(x,y)$ in terms of Legendre functions and global Poisson kernel $P^{(\mu)}(x,y)$.
\begin{theorem}
\label{thm:pk:formula}
For every $x\in D$ and $\mu>0$ we have 
\begin{equation}
\label{eq:pk:formula}
  P_D^{(\mu)}(x,y) = 
    \left\{
      \begin{array}{rl}
               \dfrac{2x_1}{(2\pi)^{n/2}}\dfrac{y_n^{\mu-\nu-1}}{x_n^{\mu+\nu}}e^{-n\pi i/2}\dfrac{Q^{\nu+1/2}_{\mu-1/2}(\cosh\rho)}{\sinh^{\nu+1/2}\rho}, & y\in\partial_1D , \\[20pt]
         \dfrac{\Gamma(\mu+\nu)}{\Gamma(\mu)\pi^\nu }\left(\dfrac{x_n^{2\mu}}{\left| x- y\right|^{2\mu+2\nu}}- \dfrac{x_n^{2\mu}}{ \left| x- \bar y\right|^{2\mu+2\nu}}\right), & y\in\partial_2D, \\

       \end{array} \right.
\end{equation}
where $\nu=\frac{n-1}{2}$, $\rho=d_{\H^n}(x,y)$.
\end{theorem}
\begin{proof}
Since $X^{(\mu)}_{\tau_D}$ takes values in  $\partial_2D$ if and only if $\tau_D=\infty$, for $A\in\mathcal{B}(\partial_2D)$ we get
\formula{\pr^x\left(X^{(\mu)}_{\tau_D}\in A\right)&=\pr^x\left(X^{(\mu)}_{\infty}\in A\right)
=\pr^x\left(\tau_D=\infty;X^{(\mu)}_\infty\in A\right)\\
&=\pr^x\left(X^{(\mu)}_\infty\in A\right)-\pr^x\left(\tau_D<\infty;X^{(\mu)}_\infty\in A\right)\/.}
Applying (\ref{eq:gpk:formula}) and the reflection principle given in Proposition \ref{prop:reflectionprinciple} we arrive at
\formula{
\pr^x\left(X^{(\mu)}_{\tau_D}\in A\right)&{=}\frac{\Gamma(\mu+\nu)}{\Gamma(\mu)\pi^\nu }\int_A\frac{x_n^{2\mu}}{ \left(\left| x- y\right|^2\right)^{\mu+\nu}}dy- \frac{\Gamma(\mu+\nu)}{\Gamma(\mu)\pi^\nu }\int_{\bar A}\frac{x_n^{2\mu}}{ \left(\left| x- y\right|^2\right)^{\mu+\nu}}dy\\
&=\int_A\frac{\Gamma(\mu+\nu)}{\Gamma(\mu)\pi^\nu }\left(\frac{x_n^{2\mu}}{\left| x- y\right|^{2\mu+2\nu}}- \frac{x_n^{2\mu}}{ \left| x- \bar y\right|^{2\mu+2\nu}}\right)dy.
}
To deal with the case $X^{(\mu)}_{\tau_D}\in \partial_1D$, (i.e. $\tau_D<\infty$) recall the representation (\ref{eq:HBM:bessel}). We denote by $\tau_1=\inf\left\{s>0:\beta_1(s)<0\right\}$ the first exit time of the one dimensional Brownian motion $\beta_1$ from the positive half-line. Since the function $t\longmapsto A_{x_n}^{(-\mu)}(t)$ is continuous, increasing and $A_{x_n}^{(-\mu)}(0)=0$ a.s., we have  $A_{x_n}^{(-\mu)}(\tau_D)=\tau_1$ whenever $\tau_D<\infty$. Thus, for every $C\in\mathcal B(\partial_1D)$, we can write
\formula{
\pr^x\left(X^{(\mu)}_{\tau_D}\in C\right)&=\pr^x\left(\tau_D<\infty,X^{(\mu)}_{\tau_D}\in C\right)\\
&=\pr\left(\tau_1<A^{(-\mu)}_{x_n}(\infty), \left(\beta_{A^{(-\mu)}_{x_n}(\tau_D)}, R^{(-\mu)}_{A^{(-\mu)}_{x_n}(\tau_D)}\right)\in C\right)\\
&=\pr\left(\tau_1<A^{(-\mu)}_{x_n}(\infty), \left(\beta_{\tau_1}, R^{(-\mu)}_{\tau_1}\right)\in C\right)\/.
}
Moreover, we have $A_{x_n}^{(-\mu)}(\infty) = T^{(-\mu)}_0$, where $T^{(-\mu)}_0=\inf\{t>0: R^{(-\mu)}_t=0\}$ is the first hitting time of zero for Bessel process $R^{(-\mu)}$, i.e. the lift-time of $R^{(-\mu)}$. Since $(\beta_{\tau_1}, R^{(-\mu)}_{\tau_1})\in C\subset\partial_1D$ implies $R^{(-\mu)}_{\tau_1}>0$ and consequently we can remove the condition $\tau_1<A^{(-\mu)}_{x_n}(\infty)$ and write
\formula{
\pr^x\left(X^{(\mu)}_{\tau_D}\in C\right) = \pr^x\left(\left(\beta_{\tau_1}, R^{(-\mu)}_{\tau_1}\right)\in C\right) = \pr^x\left(\left(\beta_1(\tau_1),...,\beta_{n-1}(\tau_1),R^{(-\mu)}(\tau_1)\right)\in C\right)
}
The density function of $\tau_1$ is given by
\formula{
\eta_{x_1}(u)=\frac{x_1}{\sqrt{2\pi}}e^{-x_1^2/2u}u^{-3/2},\quad u>0\/,
}
and since $\tau_1$ and the process $(\beta_2,\ldots,\beta_{n-1},R^{(-\mu)})$ are independent we have
\formula{
P_D^{(\mu)}(x,y)&=\int_0^\infty \frac{1}{{(2\pi u)}^{n/2-1}}\exp\left({\frac1{2u}\sum_{k=2}^{n-1}(x_k-y_k)^2}\right)g_u^{(-\mu)}(y_n)\eta_{x_1}(u)du\\
&{=}\frac{x_1x_n}{(2\pi)^{(n-1)/2}}\left(\frac{y_n}{x_n}\right)^{\mu}\int_0^\infty\exp\left(-\frac{x_ny_n}{u}\cosh \rho\right)I_{\mu}\left(\frac{x_ny_n}{u}\right)u^{-(n+3)/2}du.
}
Here the last equality follows from (\ref{eq:bessel:density}) and (\ref{eq:coshdxy}).
Substituting $u={x_ny_n}/{s}$ and then using (\ref{eq:gr6.622}) lead to
\formula{P_D^{(\mu)}(x,y) &=\frac{x_1}{2(2\pi)^{n/2-1/2}}\frac{y_n^{\mu-\nu-1}}{x_n^{\mu+\nu}}\int_0^\infty\exp\left(-s\cosh \rho\right)I_{\mu}(s)s^{(n-1)/2}ds\\
&=\frac{x_1}{(2\pi)^{n/2-1/2}}\frac{y_n^{\mu-\nu-1}}{x_n^{\mu+\nu}}\sqrt{\frac2{\pi}}e^{-n\pi i/2}\frac{Q^{n/2}_{\mu-1/2}(\cosh \rho)}{\sinh^{n/2} \rho}\\
&=\frac{2x_1}{(2\pi)^{n/2}}\frac{y_n^{\mu-\nu-1}}{x_n^{\mu+\nu}}e^{-n\pi i/2}\frac{Q^{\nu+1/2}_{\mu-1/2}(\cosh\rho)}{\sinh^{\nu+1/2}\rho}.
}
\end{proof}
The explicit formula for $P^{(\mu)}_D(x,y)$ presented above leads to the uniform estimates of the Poisson kernel $P^{(\mu)}_D(x,y)$. Note that the constant appearing in the estimates given below depends only on $n$ and the $\mu>0$. 
\begin{theorem}
\label{thm:Poisson:estimate}
There exists $c=c(n,\mu)>0$ such that for every $x\in D$, $y\in \partial D$ we have
\begin{equation*}
  P_D^{(\mu)}(x,y) \stackrel{c}{\approx} 
    \left\{
      \begin{array}{rl}
         \dfrac{x_1x_ny_n^{2\mu}}{\left|x-y\right|^{n}\left(\left|x-y\right|^2+x_ny_n\right)^{\mu+1/2}}, & y_1=0 , \\[20pt]
         \dfrac{x_1y_1x_n^{2\mu}}{\left|x-\bar y\right|^2\left(\left|x-y\right|^2\right)^{\mu+\nu}}, & y_n=0, 
       \end{array} \right.
\end{equation*}
where $\nu=\frac{n-1}{2}$ and $\rho=d_{\H^n}(x,y)$.
\end{theorem}
\begin{proof}
It is an easy consequence of formulas (\ref{eq:pk:formula}), (\ref{1/z}) and  (\ref{q}).

\end{proof}

\section{$\lambda$ - Poisson kernel}
\label{sec:LambdaPK}
In this Section we will use the notation $\ex^{(\mu)}_x$ to denote expectation with respect to the distribution of the hyperbolic Brownian motion with drift $\mu>0$ starting from $x\in \H^n$ and consequently, we will drop the superscript $(\mu)$ in the notation of the process whenever it appears under $\ex^{(\mu)}_x$. 

For any bounded domain $U \subset\subset \H^n$ we introduce the following probabilistic definition of the $\lambda$-Poisson kernel ($\lambda\geq 0$) for HBM with drift $\mu>0$
\begin{equation}
\label{eq:blpk}
P^{\mu,\lambda}_U(x,y) = \ex^{(\mu)}_x\left[e^{-\lambda \tau_U};X_{\tau_U}\in dy\right]\/,\quad x\in U\/, y\in \partial U\/.
\end{equation}
where $\tau_U$ is the first exit time of the process from the set $U$. In the following theorem we introduce the connection between $\lambda$-Poisson kernels and Poisson kernels for HBMs with different drift parameters.

\begin{theorem}
For any bounded domain $U\subset\subset \H^n$ and $\mu,\eta> 0$ we have
\label{girsanov}
\begin{equation}
\label{eq:lPoisson:bounded}
\ex^{(\mu)}_x\left[e^{-\frac12(\eta^2-\mu^2){\tau}_U};X_{\tau_U}\in dy\right]=\left(\frac{x_n}{y_n}\right)^{\mu-\eta}\ex^{(\eta)}_x\left[X_{\tau_U}\in dy\right].
\end{equation}
\end{theorem}
\begin{proof}
The last coordinate of HBM with drift is of the form $X_n^{(\mu)}(t)=x_n\exp\left(W_t-\eta t\right)$, where $W_t=B_n(t)+(\eta-\mu)t$ for one-dimensional Brownian motion $B_n$. According to Girsanov's theorem, for every $T>0$,  process $\left\{W_t, 0\leq t\leq T\right\}$ is a standard Brownian motion with respect to measure $Q_T$ defined by
\formula{
  \frac{dQ_T}{d\pr}=\exp\left((\mu-\eta)B_n(T)-\frac{(\eta-\mu)^2}2T\right)
  }
and consequently, the process $\left\{ X_t^{(\mu)}, 0\leq t \leq T\right\}$  is HBM with drift with respect to measure $Q_T$ and its generator is $\frac12\Delta_\eta$. Hence, for any $A\subset U$ we obtain
\begin{align*}
\ex^{(\eta)}_x\left[\tau_U\right.&\left.<T;X_{\tau_U}\in A \right]=\ex_x^{Q_T}\left[\tau_U<T;X_{\tau_U}\in A \right]\\
&=\ex_x^{(\mu)}\left[\tau_U<T;\exp\left((\mu-\eta)B_n\left(\tau_U\right)-\frac{(\eta-\mu)^2}2\tau_U\right);X_{\tau_U}\in A\right]\\
&=\ex_x^{(\mu)}\left[\tau_U<T;\left(\exp\left(B_n\left(\tau_U\right)-\mu\tau_U\right)\right)^{\mu-\eta}e^{-\frac12(\eta^2-\mu^2)\tau_U};X_{\tau_U}\in A\right]\\
&=x_n^{\eta-\mu}\ex_x^{(\mu)}\left[\tau_U<T;(X_n\left(\tau_U\right))^{\mu-\eta}e^{-\frac12(\eta^2-\mu^2)\tau_U};X_{\tau_U}\in A\right].
\end{align*}
Since $U$ is bounded we have $\tau_U<\infty$ a.s. and by monotone convergence theorem, taking limit as $T$ goes to infinity we obtain 
\formula{
\ex_x^{(\eta)}\left[X_{\tau_U}\in A \right]= x_n^{\eta-\mu}\ex_x^{(\mu)}\left[(X_n(\tau_U))^{\mu-\eta}e^{-\frac12(\eta^2-\mu^2)\tau_U};X_{\tau_U}\in A\right],
}
which is equivalent to (\ref{eq:lPoisson:bounded}).
\end{proof}

In the case of unbounded sets $U$ (such as hyperbolic half-space $D$) the first exit time $\tau_U$ can be infinite with positive probability and then the right-hand side of (\ref{eq:blpk}) vanishes. From the other side, taking the sequence of bounded sets $U_m\subset\subset \H^n$ with smooth boundary such that $U_m\nearrow D$ and using the formula (\ref{eq:lPoisson:bounded}) we can easily see that the $\lambda$-Poisson kernel $P^{\mu,\lambda}_{U_m}$ multiplied by $y_n^{\mu-\eta}$, where $\eta=\sqrt{2\lambda+\mu^2}$ tends to $P_D^{(\eta)}$. This observation leads us to consider the following Dirichlet problem on $D$.

For given $f\in \mathcal{C}_b(\partial D)$ and $\lambda>0$ we are looking for  $u\in\mathcal{C}^2(D)$ satisfying
	\formula[dirichlet1]{
	\left(\frac12\Delta _{\mu}u\right)(x)=\lambda u(x),\quad  x \in D\/,
	}
  such that $x_n^{\sqrt{2\lambda+\mu^2}-\mu}u(x)$ is bounded and 
		\formula[dirichlet2]{
	\lim_{\substack{x\rightarrow z\\ x\in D}}x_n^{\sqrt{2\lambda+\mu^2}-\mu}u(x)=f(z)\/,\quad z\in \partial D\/.
	}
Consequently, we extend the previously given definition of $\lambda$-Poisson kernel for bounded set $U$ into the case of unbounded hyperbolic half-space $D$.
\begin{definition}
\label{defn:lpk}
For any $\lambda\geq 0$, the $\lambda$-Poisson kernel of the set $D$ for HBM with drift $\mu>0$ is a function $P_D^{\mu, \lambda}:D\times\partial D\rightarrow[0,\infty)$ such that for given $f\in\mathcal{C}_b(\partial D)$
\formula{
u(x)=\int_{\partial D} f(y)P_D^{\mu, \lambda}(x,y)dy
}
satisfies conditions (\ref{dirichlet1}) and (\ref{dirichlet2}).
\end{definition}
Justification for the correctness of the definition as well as the formula for $P_D^{\mu, \lambda}$ are given in the following
\begin{theorem}
The function $P_D^{\mu, \lambda}$ described in Definition \ref{defn:lpk} is unique and it is given by 
\begin{align}
\label{eq:lpk:formula}
 P_D^{(\mu), \lambda}(x,y)=&\left\{
      \begin{array}{rcl}
         \dfrac{x_1}{(2\pi)^{n/2}}\dfrac{x_n^{\mu-\nu-1}}{y_n^{\alpha+\nu+1}}\,e^{-n\pi i/2}\dfrac{Q^{\nu+1/2}_{\alpha-1/2}(\cosh\rho)}{\sinh^{\nu+1/2}\rho}, & y\in\partial_1D, \\[15pt]
         \dfrac{\Gamma(\alpha+\nu)}{\Gamma(\alpha)\pi^\nu }\left(\dfrac{x_n^{\alpha+\mu}}{\left| x- y\right|^{2\alpha+2\nu}}- \dfrac{x_n^{\alpha+\mu}}{ \left| x- \bar y\right|^{2\alpha+2\nu}}\right), & y\in\partial_2D, \\
       \end{array} \right.
\end{align} where $\nu = \frac{n-1}{2}$, $\eta = \sqrt{\mu^2+2\lambda}$.
\end{theorem}
\begin{proof}
For a given $f\in\mathcal{C}_b(\partial D)$, let $u$ be a $C^2$ function on $D$ satisfying (\ref{dirichlet1}) and (\ref{dirichlet2}) and consider the sequence of bounded sets $U_m\subset\subset \H^n$ with smooth boundary such that $U_m\nearrow D$. For every $m$, there exists a unique solution $w_m\in \mathcal{C}^2(U_m)$ to the above-given Dirichlet problem with $D$ replaced by $U_m$ and $f=u\in \mathcal{C}_b(\partial U_m)$, see \cite{Folland:1992}. Moreover, it is given by 
\formula{
  w_m(x) = \ex_x^{(\mu)}\left[e^{-\lambda{\tau}_{U_m}}u\left(X_{\tau_{U_m}}\right)\right]\/.
} 
Since $u$ is also a solution, using (\ref{eq:lPoisson:bounded}) we get
\formula{
u(x) &= \ex_x^{(\mu)}\left[e^{-\lambda{\tau}_{U_m}}u\left(X_{\tau_{U_m}}\right)\right]
     = x_n^{\mu-\eta}\ex_x^{(\eta)}\left[(X_{\tau_{U_m}})^{\eta-\mu}u\left(X_{\tau_{U_m}}\right)\right]\/.
}
By dominated convergence theorem we get
\formula[eq:u:definition]{
u(x)=x_n^{\mu-\eta}\ex_x^{(\eta)}\left[f\left(X_{\tau_D}\right)\right].
} 
It gives the uniqueness of the solution of the above-given Dirichlet problem and consequently, it proves correctness of the Definition \ref{defn:lpk}. 
If we now define $u$ by (\ref{eq:u:definition}), then obviously $h(x)= x_n^{\eta-\mu}u(x) = \ex_x^{(\eta)}\left[f\left(X_{\tau_D}\right)\right]$ is bounded and the condition (\ref{dirichlet2}) holds by continuity of sample paths of the process. Since $\Delta_\eta h(x)=0$, we have
\begin{align*}
\frac12\Delta_\mu u(x)=&\ \frac12\Delta_\mu \left(x_n^{\mu-\eta}h(x)\right)\\
=&\ x_n^{\mu-\eta}\frac12\Delta_\mu h(x)+(\mu-\eta)x_n^{\mu-\eta+1}\frac{\partial h}{\partial x_n}(x)+\frac{\eta^2-\mu^2}{2}x_n^{\mu-\eta}h(x)\\
=&\ x_n^{\mu-\eta}\frac12\left(\Delta_\mu-2(\mu-\eta)x_n\frac{\partial}{\partial x_n} \right)h(x)+(\mu-\eta)x_n^{\mu-\eta+1}\frac{\partial h}{\partial x_n}(x)+\lambda x_n^{\mu-\eta}h(x)\\
=&\ x_n^{\mu-\eta}\frac12\Delta_\eta h(x)+\lambda x_n^{\mu-\eta}h(x)=\lambda u(x).
\end{align*}
Finally, the formula (\ref{eq:lpk:formula}) follows from (\ref{eq:u:definition}) and Theorem \ref{thm:pk:formula}.
\end{proof}

\section{APPENDIX}
Since many formulas appearing in the paper involve the Legendre functions $Q_b^{a}(z)$ for convenience of the Reader we collect here some basic information and properties of those functions. Here we follow the exposition given in \cite{Erdelyi:1955} (see also \cite{GradsteinRyzhik:2007}). In the second part of the section we introduce several lemmas providing some detailed properties regarding differential formulas and uniform estimates.

The \textit{Legendre functions} $Q^{a}_{b}(z)$ and $P_b^a(z)$ are solutions of the \textit{Legendre's differential equation} (\cite{GradsteinRyzhik:2007}, p.958)
\formula{
(1-z^2)u''-2zu'+(a(a+1)-b(1-z^2)^{-1})u=0,\quad b+1>a>-\frac12\/.
}
For $a>0, b\in\R$ the asymptotic behavior of the function $Q_b^a(z)$ is described by (see \cite{Erdelyi:1955}, 3.9.2 (10), (21))
\begin{eqnarray}
\label{eq:legendre:asymp01}
\lim_{z\rightarrow1^{+}}(z-1)^{a/2}Q^a_b(z)&=&e^{a \pi i}2^{a/2-1}\Gamma(a).\\
\label{eq:legendre:asymp02}
\lim_{z\rightarrow\infty}z^{b+1}Q^{a}_{b}(z)&=&\frac{\sqrt{\pi} e^{a \pi i}\Gamma(b+a+1)}{2^{b+1}\Gamma(b+3/2)}.
\end{eqnarray}
Using the above-given properties of $Q_b^a(z)$ we easily obtain 
\begin{equation}
\label{eq:legendre:asymp03}
\lim_{z\rightarrow1^{+}}(z-1)^{a}\frac{Q^a_b(z)}{(z^2-1)^{a/2}}=\frac12e^{a \pi i}\Gamma(a)\/,
\end{equation}
\begin{equation}
\label{eq:legendre:asymp04}
\lim_{z\rightarrow\infty}z^{a+b+1}\frac{Q^{a}_{b}(z)}{(z^2-1)^{a/2}}=\frac{\sqrt{\pi} e^{a \pi i}\Gamma(b+a+1)}{2^{b+1}\Gamma(b+3/2)}\/,
\end{equation}
whenever $0<a<b$. Finally, we introduce the following
\begin{lemma}
\label{lem:legendre:diff}
Let $z>1$. Then for $a,b\in \R$ we have
\formula[eq:legendre:diff]{
\frac{\partial}{\partial z}\frac{Q_b^a(z)}{(z^2-1)^{a/2}}=\frac{Q_b^{a+1}(z)}{(z^2-1)^{(a+1)/2}},
}
\end{lemma}

\begin{proof}
We apply equality (see \cite{Erdelyi:1955}, 3.8 (1) p.160)
\formula{
(b+a)(b-a-1)Q^{a-1}_b(z)=Q_b^{a+1}(z)+2a\frac{z}{\sqrt{z^2-1}}\,Q_b^{a}(z)
}
to the formula formula (see \cite{Erdelyi:1955}, 3.8 (9) p.161) 
\formula{
\frac{\partial}{\partial z}Q_b^a(z) = \frac{(a+b)(b-a+1)}{\sqrt{z^2-1}}Q^{a-1}_{b}(z) - \frac{az}{z^2-1}Q_b^a(z)
}
and we get
\formula[eq:legendre:proof01]{
\frac{\partial}{\partial z}\,Q_b^a(z)=\frac{1}{\sqrt{z^2-1}}\,Q_b^{a+1}(z)+\frac{az}{z^2-1}\,Q_b^a(z).
}
Hence
\formula{
\frac{\partial}{\partial z}\frac{Q_b^a(z)}{(z^2-1)^{a/2}} &= \frac{(z^2-1)^{a/2}\dfrac{\partial}{\partial z}Q_b^a(z)-azQ_b^a(z)(z^2-1)^{a/2-1}}{(z^2-1)^{a}}\\
&=  \frac{(z^2-1)\dfrac{\partial}{\partial z}Q_b^a(z)-azQ_b^a(z)}{(z^2-1)^{a/2+1}} = \frac{Q_b^{a+1}(z)}{(z^2-1)^{(a+1)/2}}\/,
}
where the last equality follows from (\ref{eq:legendre:proof01}).
\end{proof}

\begin{lemma}
For every $m =0,1,2,3...$ and $b\in\R$ we have
\formula[suma]{
Q_b^{m+1/2}(z) =\,i(-1)^m\sqrt{\frac{\pi}{2}}\Gamma\left(b+m+\frac 32\right) \sum_{k=0}^m\frac{a(m,k)}{\Gamma(b+3/2+k)}\frac{(z-\sqrt{z^2-1})^{b+1/2+k}}{\sqrt{z^2-1}^{1/2+k}},
}
where $a(m,k)$ are the coefficients of the Bessel polynomials 
\formula{a(m,k)=\frac{(m+k)!}{(m-k)!k!2^k}.}
\end{lemma}
\begin{proof}
To prove the result we apply mathematical induction with respect to $m$. For $m=0$ and $b\in \R$ the equality (\ref{suma}) is just the formula 3.6.1 (12) in \cite{Erdelyi:1955}.
Assume thesis for some $m\geq1$ and any $b\in\R$. 
Denoting
\formula{
f_q^p(z)=\frac{(z-\sqrt{z^2-1})^{p}}{(z^2-1)^{q}}\ \ p\in \R, q>0,
}
we have following recursion
\formula{
\frac{\partial}{\partial z}f_q^p(z)=-(p+2q)f^p_{q+1/2}(z)-2qf^{p+1}_{q+1}(z).
}
It comes from  bellow-given calculation 
\formula{
\frac{\partial}{\partial z}f_q^p(z)=&
p\frac{(z-\sqrt{z^2-1})^{p-1}}{(z^2-1)^{q}}(1-\frac{z}{\sqrt{z^2-1}})-2qz\frac{(z-\sqrt{z^2-1})^{p}}{(z^2-1)^{q+1}}\\
=&-p\frac{(z-\sqrt{z^2-1})^{p}}{(z^2-1)^{q+1/2}}-2q(z-\sqrt{z^2-1})\frac{(z-\sqrt{z^2-1})^{p}}{(z^2-1)^{q+1}}+\\
&-2q\sqrt{z^2-1}\frac{(z-\sqrt{z^2-1})^{p}}{(z^2-1)^{q+1}}\\
=&-(p+2q)\frac{(z-\sqrt{z^2-1})^{p}}{(z^2-1)^{q+1/2}}-2q\frac{(z-\sqrt{z^2-1})^{p+1}}{(z^2-1)^{q+1}}.
}
According to (\ref{eq:legendre:diff}) and inductive assumption we get
\begin{align*}
\sqrt\frac{2}{\pi}&\frac{i(-1)^{m+1}}{\Gamma\left(b+m+\frac 32\right)} \frac{Q_b^{m+1+1/2}(z)}{\sqrt{z^2-1}^{(m+1)/2+1/4}}=-\sqrt\frac{2}{\pi}\frac{i(-1)^{m}}{\Gamma\left(b+m+\frac 32\right)} \frac{\partial}{\partial z}\frac{Q_b^{m+1/2}(z)}{\sqrt{z^2-1}^{m+1/2}}\\
&=\frac{\partial}{\partial z}\left(\sum_{k=0}^m\frac{(m+k)!}{(m-k)!k!2^k\Gamma(b+3/2+k)}f^{b+\frac12+k}_{\frac{m+1}2+\frac{k}2}(z)\right)\\
&=-\sum_{k=0}^m\frac{\frac{(m+k)!}{(m-k)!k!2^k}}{\Gamma(b+\frac32+k)}\left((b+m+\frac 32+2k)f^{b+\frac12+k}_{\frac m2+1+\frac k2}(z)+(m+1+k)f^{b+\frac32+k}_{\frac m2+\frac32+\frac k2}(z)\right)\\
&=-\sum_{k=0}^m \frac{(m+k)!}{(m-k)!k!2^k\Gamma(b+\frac32+k)}(b+m+\frac 32+2k)f^{b+\frac12+k}_{\frac m2+1+\frac k2}(z)+\\
&\ \ \ -\sum_{k=1}^{m+1}\frac{(m+k-1)!}{(m-k+1)!(k-1)!2^{k-1}\Gamma(b+\frac12+k)}(m+k)f^{b+\frac12+k}_{\frac m2+1+\frac k2}(z)\\
\end{align*}
\begin{align*}
&=-\frac{b+m+\frac 32}{\Gamma(b+\frac32)}(b+m+\frac 32)f^{b+\frac12}_{\frac m2+1}(z)-\frac{(2m)!(2m+2)}{m!2^m\Gamma(b+\frac52+m)}(b+\frac32+m)f^{b+\frac32+m}_{\frac m2+1+\frac (m+1)2}(z)+\\
&\ \ \ -\sum_{k=1}^m\left(\frac{\frac{(m+k)!}{(m-k)!k!2^k}}{\Gamma(b+\frac32+k)}(b+m+\frac 32+2k)+\frac{\frac{(m+k)!}{(m-k+1)!(k-1)!2^{k-1}}}{\Gamma(b+\frac12+k)}\right)f^{b+\frac12+k}_{\frac m2+1+\frac k2}(z).\\
\end{align*}
We simplify the expression in last parenthesis as follows 
\begin{align*}
\frac{\frac{(m+k)!}{(m-k)!k!2^k}}{\Gamma(b+\frac32+k)}&(b+m+\frac 32+2k)+\frac{\frac{(m+k)!}{(m-k+1)!(k-1)!2^{k-1}}}{\Gamma(b+\frac12+k)}=\\
&=\frac{\frac{(m+k)!}{(m-k+1)!k!2^k}}{\Gamma(b+\frac32+k)}\left((b+m+\frac 32 +2k)(m-k+1)+2k(b+\frac 12+k)\right)\\
&=\frac{\frac{(m+k)!}{(m-k+1)!k!2^k}}{\Gamma(b+\frac32+k)}(b+m+\frac 32)(m+k+1)\\
&=\frac{\frac{(m+k+1)!}{(m-k+1)!k!2^k}}{\Gamma(b+\frac32+k)}(b+m+\frac 32).
\end{align*}
Consequentely
\begin{align*}
\sqrt\frac{2}{\pi}&\frac{i(-1)^{m+1}}{\Gamma\left(b+m+\frac 32\right)} \frac{Q_b^{m+1+1/2}(z)}{\sqrt{z^2-1}^{(m+1)/2+1/4}}=\\
&=-\frac{b+m+\frac 32}{\Gamma(b+\frac32)}(b+m+\frac 32)f^{b+\frac12}_{\frac m2+1}(z)-\frac{\frac{(2m)!}{m!2^m}(2m+2)}{\Gamma(b+\frac52+m)}(b+\frac32+m)f^{b+\frac32+m}_{\frac m2+1+\frac (m+1)2}(z)+\\
&\ \ \ -\sum_{k=1}^m\frac{\frac{(m+k+1)!}{(m-k+1)!k!2^k}}{\Gamma(b+\frac32+k)}(b+m+\frac 32)f^{b+\frac12+k}_{\frac m2+1+\frac k2}(z)\\
&=-(b+m+\frac 32)\sum_{k=0}^{m+1}\frac{\frac{(m+k+1)!}{(m-k+1)!k!2^k}}{\Gamma(b+\frac32+k)}f^{b+\frac12+k}_{\frac m2+1+\frac k2}(z)\\
&=-(b+m+\frac 32)\sum_{k=0}^{m+1}\frac{\frac{(m+k+1)!}{(m-k+1)!k!2^k}}{\Gamma(b+\frac32+k)}\frac{(z-\sqrt{z^2-1})^{b+\frac12+k}}{(z^2-1)^{\frac m2+1+\frac k2}}.
\end{align*}
Multiplying polar expressions by  $(z^2-1)^{m/2+3/4}\sqrt\frac{\pi}2\frac{\Gamma\left(b+m+\frac 32\right)}{i(-1)^{m+1}}$ we get formula in thesis for  $m+1$. 
\end{proof}

\begin{lemma}
\label{q-q}
For $1<x<y$ we have
\begin{equation}\label{ogr2}e^{-ai\pi }\left[ \frac{Q^a_b(x)}{(x^2-1)^{a/2}}-\frac{Q^a_b(y)}{(y^2-1)^{a/2}}\right]\stackrel{c(a,b)}{\approx}\frac{y-x}{y-1}\frac1{(x-1)^{a}x^{b+1}}.\end{equation}
\end{lemma}
 \begin{proof}
 For $x>1$ function $Q^a_b(x)$ is continuous and positive. Using equalities  (\ref{eq:legendre:asymp03}) and (\ref{eq:legendre:asymp04}) we obtain
\begin{equation}
\label{q}
e^{-ai\pi } \frac{Q^a_b(x)}{(x^2-1)^{a/2}}\stackrel{c}{\approx}\frac1{(x-1)^{a}x^{b+1}}
\end{equation}
Let assume now that  $(y-1)<2(x-1)$.  Mean value theorem and  Lemma \ref{lem:legendre:diff} gives us
$$e^{-ai\pi }\frac{\frac{Q^a_b(x)}{(x^2-1)^{a/2}}-\frac{Q^a_b(y)}{(y^2-1)^{a/2}}}{y-x}=e^{-(a+1)i\pi }\frac{Q^{a+1}_b(\xi)}{(\xi^2-1)^{(a+1)/2}},$$	
where $\xi\in[x,y]$. The right-hand side function is decreasing, hence
\begin{eqnarray*}e^{-(a+1)i\pi }\frac{Q^{a+1}_b(\xi)}{(\xi^2-1)^{(a+1)/2}}&\leq& e^{-(a+1)i\pi }\frac{Q^{a+1}_b(x)}{(x^2-1)^{(a+1)/2}}\stackrel{(\ref{q})}{\leq}\frac{c}{(x-1)^{a+1}x^{b+1}}\\
&\leq&\frac{2c}{y-1}\frac{1}{(x-1)^{a}x^{b+1}}\end{eqnarray*}
In case of $(y-1)\geq2(x-1)$  it holds that 
\begin{equation}
\label{y2x}
2\frac{y-x}{y-1}\geq1.
\end{equation}
Than
\begin{eqnarray*}e^{-ai\pi }\left[\frac{Q^a_b(x)}{(x^2-1)^{a/2}}-\frac{Q^a_b(y)}{(y^2-1)^{a/2}}\right]&\leq& e^{-ai\pi }\frac{Q^a_b(x)}{(x^2-1)^{a/2}}
\stackrel{(\ref{q})}{\leq}\frac{c}{(x-1)^{a}x^{b+1}}\\&\stackrel{(\ref{y2x})}{\leq}& \frac{2c}{y-1} \frac{y-x}{(x-1)^{a}x^{b+1}}.
\end{eqnarray*}
To prove the another inequality we estimate
\begin{align*}
e^{-ai\pi }&\left[ \frac{Q^a_b(x)}{(x^2-1)^{a/2}}-\frac{Q^a_b(y)}{(y^2-1)^{a/2}}\right]\left(\frac{y-1}{y-x}(x-1)^{a}x^{b+1}\right)\\
=\ &\frac{e^{-ai\pi }}{1-\frac{x-1}{y-1}}\left[ (x-1)^{a}x^{b+1}\frac{Q^a_b(x)}{(x^2-1)^{a/2}}-\left(\frac{x-1}{y-1}\right)^a\left(\frac{x}{y}\right)^{b+1}(y-1)^{a}y^{b+1}\frac{Q^a_b(y)}{(y^2-1)^{a/2}}\right]\\
\stackrel{(\ref q)}{\geq}&\,\frac{1}{1-\frac{x-1}{y-1}}\left[ \frac1c-c\left(\frac{x-1}{y-1}\right)^a\right]
\end{align*}
The last expression tends to $\frac1c$ as $\frac{x-1}{y-1}\rightarrow 0$. There exist $\varepsilon>0$ such that for  $\frac{x-1}{y-1}<\varepsilon$
$$e^{-ai\pi }\left[\frac{Q^a_b(x)}{(x^2-1)^{a/2}}-\frac{Q^a_b(y)}{(y^2-1)^{a/2}}\right]\left(\frac{y-1}{y-x}(x-1)^{a}x^{b+1}\right)\geq\frac1{2c},$$
what is equivalent  to
$$e^{-ai\pi }\left[\frac{Q^a_b(x)}{(x^2-1)^{a/2}}-\frac{Q^a_b(y)}{(y^2-1)^{a/2}}\right]\geq\frac1{2c}\frac{y-x}{y-1}\frac1{(x-1)^{a}x^{b+1}}.$$
For $\frac{x-1}{y-1}\geq\varepsilon$ we use mean value theorem theorem again
\begin{eqnarray*}
e^{-ai\pi }\frac{\frac{Q^a_b(x)}{(x^2-1)^{a/2}}-\frac{Q^a_b(y)}{(y^2-1)^{a/2}}}{y-x}&=&e^{-(a+1)i\pi }\frac{Q^{a+1}_b(\xi)}{(\xi^2-1)^{(a+1)/2}}\geq e^{-(a+1)i\pi }\frac{Q^{a+1}_b(y)}{(y^2-1)^{(a+1)/2}}\\
&\geq&\frac1c\frac1{(y-1)^{a+1}y^{b+1}}\geq\frac{\varepsilon^a}{c}\frac1{(y-1)(x-1)^{a}x^{b+1}}.
\end{eqnarray*}
 \end{proof}

 \begin{lemma}
For $0<x<y$, $a>0$ we have
\begin{equation}
 \label{1/z}
\frac{1}{x^a}-\frac{1}{y^a}\stackrel{c}{\approx} \frac{y-x}{yx^a}.
\end{equation}
 \end{lemma}
 \begin{proof}
Assume that  $y\geq2x$. It follows that $1\geq\frac{y-x}{y}\geq\frac12$. Hence
\begin{align*}
\frac{1}{x^a}-\frac{1}{y^a}&\leq \frac{1}{x^a}\leq2\frac{y-x}{y}\frac{1}{x^a},\\
\frac{1}{x^a}-\frac{1}{y^a}&\geq\frac{1}{x^a}-\frac{1}{(2x)^a}=\frac{2^a-1}{2^a}\frac{1}{x^a}\geq\frac{2^a-1}{2^a}\frac{y-x}{y}\frac{1}{x^a}.
\end{align*}
If  $y<2x$ we use mean value theorem
 $$\frac{\frac{1}{x^a}-\frac{1}{y^a}}{y-x}=a\frac1{z^{a+1}},$$
 where $z\in[x,y]$. Then we get
\begin{align*}
 \frac1{z^{a+1}}&\leq\frac1{x^{a+1}}\leq\frac2{yx^a},\\
 \frac1{z^{a+1}}&\geq\frac1{y^{a+1}}\geq\frac1{2^a}\frac1{yx^a}.
\end{align*}
 
 \end{proof}

\bibliography{bibliography}
\bibliographystyle{plain}

\end{document}